# MATHEMATICAL PROPERTIES AND FINITE-POPULATION CORRECTION FOR THE WILSON SCORE INTERVAL

B. O'Neill,[*] *Australian National University*[**]

WRITTEN 28 AUGUST 2021

Abstract

In this paper we examine the properties of the Wilson score interval, used for inferences for an unknown binomial proportion parameter. We examine monotonicity and consistency properties of the interval and we generalise it to give two alternative forms for inferences undertaken in a finite population. We discuss the nature of the "finite population correction" in these generalised intervals and examine their monotonicity and consistency properties. This analysis gives the appropriate confidence interval for an unknown population proportion or unknown unsampled proportion in a finite or infinite population. We implement the generalised confidence interval forms in a user-friendly function in **R**.

BINOMIAL PROPORTION; WILSON SCORE INTERVAL; MONOTONICITY; CONSISTENCY; FINITE POPULATION CORRECTION.

Interval estimation for an unknown proportion parameter using binary data is a well-known statistical problem. It is commonly covered in introductory statistical courses and is considered to be an essential tool in the applied statistician's toolkit. Wilson (1927) introduced a form of confidence interval for this problem that has now become known by the eponymous name of the "Wilson score interval". This is a simple confidence interval that is one of the key tools used in this kind of inference problem. The performance of this interval has been examined and contrasted with other confidence intervals for a population proportion and it is known to perform well in terms of accuracy and coverage in a wide class of cases (see e.g., Ghosh 1979; Blyth and Still 1983; Vollset 1993; Agresti and Coull 1998; Brown, Cai and DasGupta 2001).

## 1. The Wilson score interval and its monotonicity and consistency properties

The Wilson score interval arises in the context of the binomial model, best characterised by an underlying sequence of exchangeable binary values. Consider a population composed of IID binary values $X_1, \ldots, X_N \sim \text{IID Bern}(\theta)$ with population size $N \in \overline{\mathbb{N}} \equiv \{1, 2, 3, \ldots, \infty\}$ (i.e., we allow either a finite or infinite population). We take a sample of size $1 \leq n \leq N$ from this population, and without loss of generality we take these to be the values $X_1, \ldots, X_n$. To facilitate our analysis of this general case (allowing for a finite or infinite population), we will make use of the sample proportion $\bar{X}_n$, the population proportion $\bar{X}_N$ and the unsampled proportion $\bar{X}_{n:N}$, which are defined respectively by:

$$\bar{X}_n \equiv \frac{1}{n}\sum_{i=1}^{n} X_i \qquad \bar{X}_N \equiv \frac{1}{N}\sum_{i=1}^{N} X_i \qquad \bar{X}_{n:N} \equiv \frac{1}{N-n}\sum_{i=n+1}^{N} X_i.$$

---

[*] E-mail address: ben.oneill@hotmail.com.
[**] Research School of Population Health, Australian National University, Canberra ACT 0200, Australia.

From the law-of-large-numbers we can write the superpopulation parameter $\theta$ in terms of limits of the latter quantities as $\theta = \lim_{N\to\infty} \bar{X}_N = \lim_{N\to\infty} \bar{X}_{n:N}$.[1] The **Wilson score interval** for $\theta$ (Wilson 1927) at confidence level $1-\alpha$ can be written as the following function:

$$\text{CI}_\infty \equiv \text{CI}_\infty(\alpha, n, \bar{x}) \equiv \left[\frac{n\bar{x} + \tfrac{1}{2}\chi_\alpha^2}{n + \chi_\alpha^2} \pm \frac{\chi_\alpha}{n + \chi_\alpha^2}\sqrt{n\bar{x}(1-\bar{x}) + \tfrac{1}{4}\chi_\alpha^2}\right],$$

where $\chi_\alpha^2$ is the critical point of the chi-squared distribution with one degree-of-freedom (with upper tail area $\alpha$). In this formulate of the function we have treated the sample size $n$ and the sample proportion $\bar{x}$ separately to facilitate general analysis. Substituting the observed values of the sample size and sample proportion, and the chosen value of $\alpha$ for the confidence level, gives the confidence interval. (Our notation using the subscript infinity is chosen to indicate that this is the confidence interval for inference about the parameter $\theta$, which is the population proportion for an infinite population.)

In this paper we will generalise this well-known confidence interval to obtain intervals for the quantities $\bar{X}_N$ and $\bar{X}_{n:N}$ (both of which reduce down to the parameter $\theta$ in the limit as $N \to \infty$). This will give us forms for the "finite-population correction" for the interval, allowing use for inferences of the population proportion or unsampled proportion from a finite or infinite population. Before attempting to generalise the Wilson score interval for a finite population we will first examine some basic intuitive properties of this interval that we would like to preserve in our generalisation.

The standard Wilson score interval shown above has a number of sensible monotonicity and consistency properties which satisfy basic intuitive requirements for an inference about $\theta$. To examine these, we will first define some useful functions pertaining to the confidence interval. We define the **lower and upper bounds** of the interval respectively as:

---

[1] Here we use the Cesàro limits, but these do not always exist so the parameter $\theta$ can be defined as the Banach limit of the observable sequence (see O'Neill 2009 for further discussion on this matter).

$$L_\infty(\alpha, n, \bar{x}) \equiv \frac{n\bar{x} + \frac{1}{2}\chi_\alpha^2}{n + \chi_\alpha^2} - \frac{\chi_\alpha}{n + \chi_\alpha^2}\sqrt{n\bar{x}(1 - \bar{x}) + \frac{1}{4}\chi_\alpha^2},$$

$$U_\infty(\alpha, n, \bar{x}) \equiv \frac{n\bar{x} + \frac{1}{2}\chi_\alpha^2}{n + \chi_\alpha^2} + \frac{\chi_\alpha}{n + \chi_\alpha^2}\sqrt{n\bar{x}(1 - \bar{x}) + \frac{1}{4}\chi_\alpha^2},$$

and we define the **width function** for the interval as:

$$w_\infty(\alpha, n, \bar{x}) \equiv U_\infty(\alpha, n, \bar{x}) - L_\infty(\alpha, n, \bar{x}) = \frac{2\chi_\alpha}{n + \chi_\alpha^2}\sqrt{n\bar{x}(1 - \bar{x}) + \frac{1}{4}\chi_\alpha^2}.$$

In practice, the sample size $n$ is a non-negative integer and the sample proportion $\bar{x}$ is restricted to the possible values $0, 1/n, \ldots, (n-1)/n, 1$. Moreover, we must have $n \geq 1$ (i.e., at least one data point) in order to get any non-vacuous inference. Notwithstanding these constraints, for our purposes it is more useful to define these functions for all real values $n \geq 0$ and all real values $0 \leq \bar{x} \leq 1$ so that we can examine the functions over a broader range than the values that appear in practical cases. (By using real values we also allow ourselves to examine the properties of these functions through derivatives.) It is also notable that we treat the sample proportion as a separate input from the sample size in these functions, even though they are related in practice (i.e., when we change the sample size we change the sample, so the sample proportion can change). Defining the above functions in this way is useful because it allows us to examine properties of the confidence interval as we change one of the inputs, holding all other inputs constant.

The main intuitive properties of the Wilson score interval are shown in Theorems 1-4 below. From the limiting results in Theorem 1 we can see that the interval is vacuous when $\alpha = 0$ (i.e., a confidence level of 100%) or $n = 0$ (i.e., with no data), which accords with our intuition. Similarly, we can see that the interval converges to a single point when $\alpha = 1$ (i.e., a confidence level of 0%) or $n \to \infty$ (i.e., we observe the entire superpopulation), which also accords with our intuition. From the monotonicity results in Theorem 2 we can see that the interval becomes more accurate as we decrease the confidence level or as we get more data, both of which accord with intuition.[2] From the monotonicity results in Theorem 3 we can see that the bounds of the interval move up or down as the sample proportion moves up or down, which reflects the fact that a higher sample proportion is *ipso facto* evidence of a higher superpopulation proportion. This is closely related to the monotone-likelihood-ratio property of the binomial distribution (Karlin and Rubin 1956). The concavity result in Theorem 3 shows that the width function is strictly concave — in fact, the symmetrical nature of the width function means that it is maximised when $\bar{x} = \frac{1}{2}$ and minimised when $\bar{x} = 0$ or $\bar{x} = 1$. In Theorem 4 we confirm that the Wilson score interval is always concentrated on the possible range of values of $\theta$ (unlike some other intervals — e.g., the Wald interval). Finally, in Theorem 5 we show explicit bounds on the width for any fixed values of $\alpha$ and $n$. Since the width is symmetric and strictly concave with respect to $\bar{x}$, the minimum width occurs at the extremes $\bar{x} = 0$ and $\bar{x} = 1$ and the maximum width occurs at the mid-point $\bar{x} = \frac{1}{2}$. This result allows us to narrow down the accuracy of the interval *a priori* for the purposes of sample size calculations.

---

[2] Some readers may balk at the idea that interval should become more accurate as we decrease the confidence level (or that this result corresponds to the theorem). To see that this corresponds to the theorem, observe that as $\alpha$ increases the confidence level *decreases*, and the width of the interval decreases (i.e., the interval becomes *more* accurate). The intuition for this is that, for a fixed set of data, there is a trade-off between confidence and accuracy in the interval — we can either have a more accurate interval with less confidence, or a less accurate interval with more confidence.

**THEOREM 1:** We have:

$$\lim_{\alpha \to 0} \text{CI}_\infty(\alpha, n, \bar{x}) = [0,1] \qquad \lim_{\alpha \to 1} \text{CI}_\infty(\alpha, n, \bar{x}) = [\bar{x}],$$

$$\lim_{n \to 0} \text{CI}_\infty(\alpha, n, \bar{x}) = [0,1] \qquad \lim_{n \to \infty} \text{CI}_\infty(\alpha, n, \bar{x}) = [\bar{x}].$$

**THEOREM 2:** For all $0 < \alpha < 1$ and $n > 0$ we have:

$$\frac{\partial w_\infty}{\partial \alpha}(\alpha, n, \bar{x}) < 0 \qquad \frac{\partial w_\infty}{\partial n}(\alpha, n, \bar{x}) < 0.$$

**THEOREM 3:** For all $0 < \alpha < 1$ and $n > 0$ we have:

$$\frac{\partial L_\infty}{\partial \bar{x}}(\alpha, n, \bar{x}) > 0 \qquad \frac{\partial U_\infty}{\partial \bar{x}}(\alpha, n, \bar{x}) > 0 \qquad \frac{\partial^2 w_\infty}{\partial \bar{x}^2}(\alpha, n, \bar{x}) < 0.$$

**THEOREM 4:** For all values of $0 \leq \alpha \leq 1$, $n \geq 0$ and $0 \leq \bar{x} \leq 1$ we have:

$$0 \leq L_\infty(\alpha, n, \bar{x}) \leq U_\infty(\alpha, n, \bar{x}) \leq 1.$$

**THEOREM 5:** For all values of $0 \leq \alpha \leq 1$ and $n \geq 0$ we have:

$$\frac{\chi_\alpha^2}{n + \chi_\alpha^2} \leq w_\infty(\alpha, n, \bar{x}) \leq \frac{\chi_\alpha}{\sqrt{n + \chi_\alpha^2}}.$$

The above statistical properties give us confidence that the Wilson score interval operates in a sensible way with respect to the confidence level, sample size, and sample proportion. The fact that the interval is concentrated on the possible values of the parameter for all possible inputs is particularly desirable. The performance of the Wilson score interval is examined in detail in Brown, Cai and DasGupta (2001) and is contrasted with other confidence intervals for an unknown proportion parameter. This analysis shows that the interval performs well in comparison to other commonly used intervals; in particular, it outperforms the Wald interval in terms of accuracy and coverage in most cases. Although the Wilson score interval is derived from a pivotal quantity with a distribution that relies on the central limit theorem (i.e., assuming $n$ is large), we can see that the confidence interval obeys sensible properties even for small values of $n$, and in fact, the interval performs surprisingly well even for quite small values. Some relevant discussion on asymptotic performance and possible improvements is found in Hall (1982), Brown, Cai, and DasGupta (2002) and Brown, Cai, and DasGupta (2003).

In addition to demonstrating desirable statistical properties, another interesting aspect of the monotonicity and concavity results in Theorem 3 is that they establish that the bounds on the width of the interval (taking $\alpha$ and $n$ as fixed) are given at the extremes $\bar{x} = 0$ and $\bar{x} = 1$. At these respective extremes the Wilson score interval reduces to:

$$\text{CI}_\infty(\alpha, n, 0) = \left[0, \frac{\chi_\alpha^2}{n + \chi_\alpha^2}\right],$$

$$\text{CI}_\infty(\alpha, n, 1) = \left[\frac{n}{n + \chi_\alpha^2}, 1\right].$$

It is notable that these special cases bear a resemblance to the heuristic "rule of three" for 95% confidence intervals in which no positive/negative outcomes are observed (see e.g., Jovanovic and Levy 1997; Tuyl, Gerlach and Mengersen 2009). Taking $\alpha = 0.05$ for a 95% confidence interval reduces the Wilson score interval to:

$$\text{CI}_\infty(\alpha, n, 0) = \left[0, \frac{3.841459}{n + 3.841459}\right],$$

$$\text{CI}_\infty(\alpha, n, 1) = \left[\frac{n}{n + 3.841459}, 1\right].$$

This form bears a clear resemblance to the corresponding intervals $[0, 3/n]$ and $[1 - 3/n, 1]$ using the "rule of three". The superiority of the Wilson score interval can be seen here from the fact that it is applicable for values of $n < 3$, unlike the latter interval.

## 2. Finite-population correction for the Wilson score interval

In practical sampling problems the population of interest is often finite, and in this case we will generally prefer to make an inference about the population proportion $\bar{X}_N$ or the unsampled proportion $\bar{X}_{n:N}$ instead of the (hypothesised) parameter $\theta$. A confidence interval for a finite population proportion is desirable when we seek knowledge about the actual finite population of interest, and a confidence interval for the unsampled proportion from a finite population is desirable when we want to exclude the sampled items from the scope of analysis. The latter can occur in cases where the process of measurement in an experiment destroys the sampled items, so that the unsampled items form the post-experiment population.[3] It can also occur when we have other reasons to exclude the sampled items from consideration.

In order to determine the "finite-population correction" for the Wilson score interval we will consider a population of size $N \geq 1$ and derive intervals for $\bar{X}_N$ and $\bar{X}_{n:N}$ from scratch. Our analysis here relies on moment results examined in O'Neill (2014). That paper establishes some basic moment results and pivotal quantities involving the above quantities for a finite population. To facilitate our analysis we will define the **effective sample sizes**:

$$n_* \equiv n \cdot \frac{N-1}{N-n} \qquad n_{**} \equiv n \cdot \frac{N-n}{N-1}.$$

The first of these will be shown to be the effective sample size for an inference about the population proportion $\bar{X}_N$ and the second will be shown to be the effective sample size for an inference about the unsampled proportion $\bar{X}_{n:N}$. If $n$ and $N-n$ are large we can invoke the central limit theorem to approximate the binomial by the normal, and using Results 14-15 in O'Neill (2014) (pp. 285-286) with slight modification we obtain the pivotal quantities:

$$\frac{n_*(\bar{X}_n - \bar{X}_N)^2}{\bar{X}_N(1-\bar{X}_N)} \overset{\text{Approx}}{\sim} \text{ChiSq}(1),$$

$$\frac{n_{**}(\bar{X}_n - \bar{X}_{n:N})^2}{\bar{X}_{n:N}(1-\bar{X}_{n:N})} \overset{\text{Approx}}{\sim} \text{ChiSq}(1).$$

We will use the above pivotal quantities to obtain versions of the Wilson score interval for the population proportion and the unsampled proportion.

---

[3] A simple example from an industrial context occurs when quality-control inspectors randomly sample cans of food and inspect the food for quality (pass/fail). In this case, opening the can in order to inspect the food means that the item is no longer saleable, and so the sampled items no longer form part of the saleable population of goods. In this type of problem, an inference about the unsampled proportion is desirable.

We will illustrate this by showing the derivation of the confidence interval for the population proportion $\bar{X}_N$. Let $\chi_\alpha^2$ denote the critical point of the chi-squared distribution with one degree-of-freedom (with upper tail area $\alpha$) we then write:

$$\begin{aligned}
1 - \alpha &\approx \mathbb{P}(n_*(\bar{X}_n - \bar{X}_N)^2 \leq \bar{X}_N(1 - \bar{X}_N)\chi_\alpha^2) \\
&= \mathbb{P}(n_*(\bar{X}_n^2 - 2\bar{X}_n\bar{X}_N + \bar{X}_N^2) \leq (\bar{X}_N - \bar{X}_N^2)\chi_\alpha^2) \\
&= \mathbb{P}\big((n_* + \chi_\alpha^2)\bar{X}_N^2 - (2n_*\bar{X}_n + \chi_\alpha^2)\bar{X}_N + n_*\bar{X}_n^2 \leq 0\big) \\
&= \mathbb{P}\left(\bar{X}_N^2 - \frac{2n_*\bar{X}_n + \chi_\alpha^2}{n_* + \chi_\alpha^2} \cdot \bar{X}_N + \frac{n_*\bar{X}_n^2}{n_* + \chi_\alpha^2} \leq 0\right) \\
&= \mathbb{P}\left(\left(\bar{X}_N - \frac{n_*\bar{X}_n + \tfrac{1}{2}\chi_\alpha^2}{n_* + \chi_\alpha^2}\right)^2 \leq \frac{\chi_\alpha^2}{n_* + \chi_\alpha^2} \cdot \frac{n_*\bar{X}_n(1 - \bar{X}_n) + \tfrac{1}{4}\chi_\alpha^2}{n_* + \chi_\alpha^2}\right) \\
&= \mathbb{P}\left(\bar{X}_N \in \left[\frac{n_*\bar{X}_n + \tfrac{1}{2}\chi_\alpha^2}{n_* + \chi_\alpha^2} \pm \frac{\chi_\alpha}{n_* + \chi_\alpha^2}\sqrt{n_*\bar{X}_n(1 - \bar{X}_n) + \tfrac{1}{4}\chi_\alpha^2}\right]\right).
\end{aligned}$$

Substituting the observed data $\bar{x}$ gives the confidence interval:

$$\mathrm{CI}_N = \mathrm{CI}_N(\alpha, n, \bar{x}) = \left[\frac{n_*\bar{x} + \tfrac{1}{2}\chi_\alpha^2}{n_* + \chi_\alpha^2} \pm \frac{\chi_\alpha}{n_* + \chi_\alpha^2}\sqrt{n_*\bar{x}(1 - \bar{x}) + \tfrac{1}{4}\chi_\alpha^2}\right].$$

From this form we can see that we have the simple relationship $\mathrm{CI}_N(\alpha, n, \bar{x}) = \mathrm{CI}_\infty(\alpha, n_*, \bar{x})$ relating the Wilson score interval for a finite population to the standard Wilson score interval. In the case where we have an infinite population we take $N = \infty$ which gives $n_* = n$, so the generalised interval reduces down to the standard form in this case. In the case of a full census of the population we take $N = n$ which gives $n_* = \infty$, so the confidence interval reduces down to $\mathrm{CI}_n(\alpha, n, \bar{x}) = [\bar{x}]$ which is just what we expect when the entire population is observed.

The confidence interval for $\bar{X}_{n:N}$ is analogous, but we use the effective sample size $n_{**}$ instead of $n_*$. Using the same working this latter interval is:

$$\mathrm{CI}_{n:N} = \mathrm{CI}_{n:N}(\alpha, n, \bar{x}) = \left[\frac{n_{**}\bar{x} + \tfrac{1}{2}\chi_\alpha^2}{n_{**} + \chi_\alpha^2} \pm \frac{\chi_\alpha}{n_{**} + \chi_\alpha^2}\sqrt{n_{**}\bar{x}(1 - \bar{x}) + \tfrac{1}{4}\chi_\alpha^2}\right].$$

From this form we can see that we have the simple relationship $\mathrm{CI}_{n:N}(\alpha, n, \bar{x}) = \mathrm{CI}_\infty(\alpha, n_{**}, \bar{x})$ relating the Wilson score interval for a unsampled part of the population to the standard Wilson score interval. In the case where we have an infinite superpopulation we take $N = \infty$ which gives $n_{**} = n$, so the generalised interval reduces down to the standard form in this case. In the case where we have a full census of the population we take $N = n$ which gives $n_{**} = 0$, so the confidence interval reduces down to $\mathrm{CI}_{n:n}(\alpha, n, \bar{x}) = [0, 1]$, reflecting the fact that there is no unsampled part of the population in this case.

An alternative succinct form for the above intervals —where all the information other than the sample proportion is incorporated into a single parameter— is provided in Appendix II. Both of the generalised intervals are implemented in the **CONF.prop** function in **stat.extend** package in **R** (O'Neill and Fultz 2020).[4] An example of use of the function is shown in the code below (inputs are shown in black and outputs are shown in blue).

```
#Load the package
library(stat.extend)

#Set the parameters
theta <- 0.7
N     <- 200
n     <- 60

#Create binary data (population, sample and unsampled)
set.seed(1956279346)
POPULATION  <- sample(x = c(0,1), size = N, replace = TRUE,
                      prob = c(1-theta, theta))
SAMPLE      <- POPULATION[1:n]
UNSAMPLED   <- POPULATION[(n+1):N]

#Compute true proportions in the population, sample, and unsampled parts
(PROP.SAMP   <- mean(SAMPLE))
[1] 0.65
(PROP.UNSAMP <- mean(UNSAMPLED))
[1] 0.7285714
(PROP.POP    <- mean(POPULATION))
[1] 0.705

#Compute 95% confidence interval for population proportion
CONF.prop(alpha = 0.05, x = SAMPLE, N = N, unsampled = FALSE)

        Confidence Interval (CI)

95.00% CI for proportion for population of size 200
Interval uses 60 binary data points from data SAMPLE with sample
proportion = 0.6500

[0.544301577788208, 0.742768160810678]

#Compute 95% confidence interval for population proportion
CONF.prop(alpha = 0.05, x = SAMPLE, N = N, unsampled = TRUE)

        Confidence Interval (CI)

95.00% CI for proportion for unsampled population of size 140
Interval uses 60 binary data points from data SAMPLE with sample
proportion = 0.6500

[0.49916421640973, 0.775811359434426]
```

---

[4] This package also contains various other confidence interval functions that accommodate finite populations and inferences for the full population or the unsampled part.

In the above example the two confidence intervals computed for the population proportion and unsampled proportion both include the true values of the quantities of interest. As can be seen, the function operates in a simple user-friendly manner where the user specified the population size `N` and a logical value `unsampled`. (The binary sample data can be input directly as the input `x` or the user can instead use the inputs `n` and `sample.prop` to specify the sample size and sample proportion in the data; in the above code we have input the full sample data.) The output of the confidence interval function gives a closed interval in the form used in the `sets` package (Meyer, Hornik and Buchta 2017), with some additional printed information about the quantity of interest in the inference.

### 3. Properties of the generalised Wilson score intervals

Since the finite-population correct for the Wilson score interval only entails a change of the effective sample size (to another value that is a non-negative real number) these generalised versions of the interval preserve all the basic coherence properties that apply to the standard version of the Wilson score interval. The relationship with respect to the sample size $n$ is now changed through the fact that this input is mediated through the effective sample size function. In Theorem 6 below we shown the general monotonicity, directional and limiting properties of the two effective sample size functions.

**THEOREM 6:** For all $0 < n < N$ and $N > 1$ we have:[5]

$$\frac{\partial n_*}{\partial n}(n,N) > 0 \qquad \frac{\partial n_*}{\partial N}(n,N) < 0,$$

$$\frac{\partial n_{**}}{\partial n}(n,N) \sim \tfrac{1}{2}N - n \qquad \frac{\partial n_{**}}{\partial N}(n,N) > 0,$$

For all $0 \leq n \leq N$ and $N \geq 1$ we have:

$$n_*(0,N) = 0 \qquad n_*(N,N) = \infty,$$
$$n_*(n,n) = \infty \qquad n_*(n,\infty) = n,$$
$$n_{**}(0,N) = 0 \qquad n_{**}(N,N) = 0,$$
$$n_{**}(n,n) = 0 \qquad n_{**}(n,\infty) = n.$$

---

[5] In these theorems we use the ~ sign to represent sign-equivalence (i.e., $a \sim b$ means sgn $a$ = sgn $b$). This is a valid equivalence relation and is used here for brevity for the directional results for the derivatives.

As with the standard Wilson score interval we can define lower and upper bound functions and the width function for each of the generalised intervals. In Theorems 7-10 we then show the montonicity and consistency properties for the generalised versions of the Wilson score interval (note that our notation for the lower and upper bound functions and width function is just as before except that we use different subscripts corresponding to the above confidence intervals). We can see from these theorems that the interval for a finite population proportion is consistent, as we would expect. This interval has the same essential properties as the standard Wilson score interval. From Theorem 8 we can see that the confidence interval becomes less accurate as the population size increases (holding the sample size fixed). This reflects the fact that with a smaller population the sample takes up a larger proportion, and since the latter is known, it gives direct information about the population proportion.

**THEOREM 7:** We have:

$$\lim_{\alpha \to 0} \text{CI}_N(\alpha, n, \bar{x}) = [0,1] \qquad \lim_{\alpha \to 1} \text{CI}_N(\alpha, n, \bar{x}) = [\bar{x}],$$

$$\lim_{n \to 0} \text{CI}_N(\alpha, n, \bar{x}) = [0,1] \qquad \lim_{n \to N} \text{CI}_N(\alpha, n, \bar{x}) = [\bar{x}],$$

$$\lim_{\alpha \to 0} \text{CI}_{n:N}(\alpha, n, \bar{x}) = [0,1] \qquad \lim_{\alpha \to 1} \text{CI}_{n:N}(\alpha, n, \bar{x}) = [\bar{x}],$$

$$\lim_{n \to 0} \text{CI}_{n:N}(\alpha, n, \bar{x}) = [0,1] \qquad \lim_{n \to N} \text{CI}_{n:N}(\alpha, n, \bar{x}) = [0,1].$$

**THEOREM 8:** For all $0 < \alpha < 1$, $0 < n < N$ and $N > 1$ we have:

$$\frac{\partial w_N}{\partial \alpha}(\alpha, n, \bar{x}) < 0 \qquad \frac{\partial w_N}{\partial n}(\alpha, n, \bar{x}) < 0 \qquad \frac{\partial w_N}{\partial N}(\alpha, n, \bar{x}) > 0.$$

$$\frac{\partial w_{n:N}}{\partial \alpha}(\alpha, n, \bar{x}) < 0 \qquad \frac{\partial w_{n:N}}{\partial n}(\alpha, n, \bar{x}) \sim n - \tfrac{1}{2}N \qquad \frac{\partial w_{n:N}}{\partial N}(\alpha, n, \bar{x}) < 0.$$

**THEOREM 9:** For all $0 < \alpha < 1$, $0 < n < N$ and $N > 1$ we have:

$$\frac{\partial L_N}{\partial \bar{x}}(\alpha, n, \bar{x}) > 0 \qquad \frac{\partial U_N}{\partial \bar{x}}(\alpha, n, \bar{x}) > 0 \qquad \frac{\partial^2 w_N}{\partial \bar{x}^2}(\alpha, n, \bar{x}) < 0,$$

$$\frac{\partial L_{n:N}}{\partial \bar{x}}(\alpha, n, \bar{x}) > 0 \qquad \frac{\partial U_{n:N}}{\partial \bar{x}}(\alpha, n, \bar{x}) > 0 \qquad \frac{\partial^2 w_{n:N}}{\partial \bar{x}^2}(\alpha, n, \bar{x}) < 0.$$

**THEOREM 10:** For all values of $0 \leq \alpha \leq 1$, $0 \leq n \leq N$, $N \geq 1$ and $0 \leq \bar{x} \leq 1$ we have:

$$0 \leq L_N(\alpha, n, \bar{x}) \leq U_N(\alpha, n, \bar{x}) \leq 1,$$
$$0 \leq L_{n:N}(\alpha, n, \bar{x}) \leq U_{n:N}(\alpha, n, \bar{x}) \leq 1.$$

**THEOREM 11:** For all values of $0 \leq \alpha \leq 1$, $0 \leq n \leq N$, $N \geq 1$ and $0 \leq \bar{x} \leq 1$ we have:

$$\frac{\chi_\alpha^2}{\tfrac{1}{4}N + \chi_\alpha^2} \leq w_{n:N}(\alpha, n, \bar{x}) \leq 1.$$

The confidence interval $\text{CI}_{n:N}$ is somewhat unusual, insofar as it does not reduce to a point-mass estimator when $n = N$ (i.e., when we take a full census of the population). This is due to an interesting trade-off that occurs when making an inference about the unsampled proportion of the population. As we put more data in the sample, we obtain more statistical information about the unsampled part of the population, but we also reduce the corresponding size of the unsampled portion of the population and thereby increase the variability of the unsampled proportion of the population. If we take too much data, the increase in the variability of the unsampled proportion outweighs the statistical information gained from the additional data and our inference gets worse. This trade-off leads to a limitation in the accuracy of the confidence interval for any fixed population size and confidence level. Specifically, from Theorem 11 we see that the width of the confidence interval for the unsampled proportion cannot be reduced below the lower bound $\chi_\alpha^2/(\tfrac{1}{4}N + \chi_\alpha^2)$, irrespective of the sample size and observed sample proportion. This lower bound is decreasing in $N$, so there are strict limitations in the accuracy of an inference about the unsampled proportion in small populations.

## 4. Isoquants and sample size calculations

In order to get a good sense of the information properties for an inference about the population mean, it is useful to examine the "isoquants" that determine the required sample size needed to give an inference with some fixed level of confidence and accuracy, as a function of the size of the size of the unsampled group. This is particularly easy in the context of the generalised Wilson score interval because the sample size enters the interval formula only through the effective sample size functions. With sample size $n$ and unsampled size $m \geq 1$ we have population size $N = n + m \geq 1$, which gives:

$$n_* = \frac{n(n + m - 1)}{m}.$$

Holding the effective sample size for the interval fixed, the isoquants are given by the function:

$$n = \frac{\sqrt{m^2 + (4n_* - 2)m + 1}}{2} - \frac{m-1}{2}.$$

In Figure 1 below we plot the isoquants for various values for the effective sample size. It can easily be shown that the function is strictly increasing and concave. This reflects that fact that increases in the size of the unsampled population require marginally less values in the sample to give an inference with fixed accuracy and confidence. (The corresponding form for the isoquants for inference about the unsampled mean is also quite simple, though less intuitive; due to limitations on the effective sample size for that inference, the latter function is defined only for values $m > n_{**}$, and for these values the function strictly decreasing and convex.)

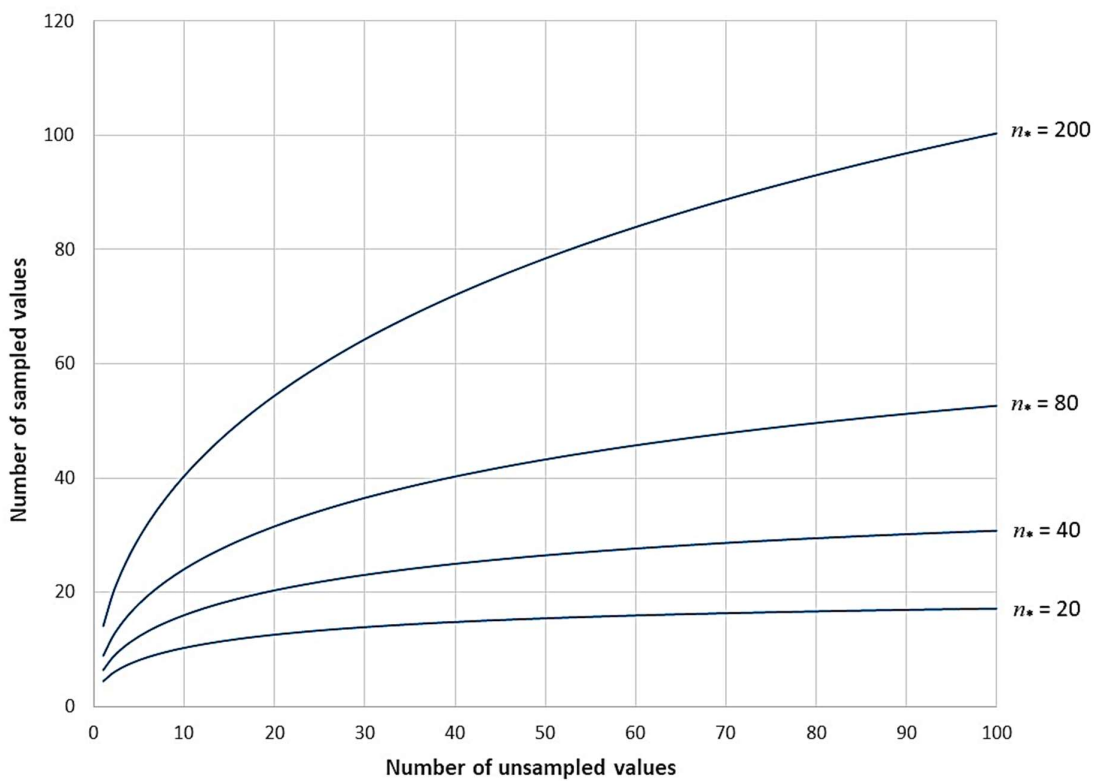

**Figure 1:** Isoquant function for inference for the population mean

An interesting aspect of sampling problems —well known to statistical educators— is that it is common for the layperson to believe that what matters to the accuracy of an inference about a population characteristic is the proportion of values in the population that are included in the sample. This lay view looks only at the direct information about the population, and ignores the statistical information the sample gives about unsampled values. It is instructive to note that if the lay understanding were true, the isoquants in the figure would be linear functions emanating from the origin and heading upward at a fixed sampling proportion. The strictly concave nature of the isoquants falsifies this view and reflects the statistical information from the sample values to the unsampled values.

The width functions for the confidence intervals in the previous sections measure the accuracy of the inference, and we have seen that these functions have various intuitive properties. It is useful to examine the required sample size for an inference with a particular level of accuracy. For any real inputs $0 < w < 1$, $0 < \alpha < 1$ and $0 \leq \bar{x} \leq 1$ we define the **required sample size** function by $\hat{n} \equiv \hat{n}(w, \alpha, \bar{x}) \equiv \inf\{n \geq 0 | w_\infty(\alpha, n, \bar{x}) \leq w\}$. It is simple to write this function in explicit form; setting $w_\infty(\alpha, n, \bar{x}) \leq w$ and simplifying for $n$ gives:

$$n^2 + 2\chi_\alpha^2 \cdot \frac{w^2 - 2\bar{x}(1-\bar{x})}{w^2} \cdot n - \chi_\alpha^4 \cdot \frac{1 - w^2}{w^2} \geq 0.$$

This quadratic inequality simplifies to:

$$\left| \frac{n}{\chi_\alpha^2} + \frac{w^2 - 2\bar{x}(1-\bar{x})}{w^2} \right| \geq \frac{\sqrt{w^2 - 4\bar{x}(1-\bar{x})w^2 + 4\bar{x}^2(1-\bar{x})^2}}{w^2},$$

and taking the infimum of all $n$ satisfying this inequality then yields the explicit solution:

$$\hat{n}(w, \alpha, \bar{x}) = \frac{\chi_\alpha^2}{w^2} \left[ \sqrt{w^2 - 4\bar{x}(1-\bar{x})w^2 + 4\bar{x}^2(1-\bar{x})^2} - w^2 + 2\bar{x}(1-\bar{x}) \right]$$

$$= \frac{\chi_\alpha^2}{w^2} \left[ \begin{array}{c} \sqrt{w^2 - 4\bar{x}(1-\bar{x})w^2 + 4\bar{x}^2(1-\bar{x})^2} \\ -\sqrt{w^4 - 4\bar{x}(1-\bar{x})w^2 + 4\bar{x}^2(1-\bar{x})^2} \end{array} \right].$$

Of course, in practice, sample size calculations must occur before observation of the sample, so the above function merely gives a range of values for different possible values of the sample proportion $\bar{x}$. The range of possible values taken over this variable is shown in Theorem 12. The upper bound in the theorem gives a conservative form for the required sample size, which is sufficient to obtain the desired level of accuracy irrespective of the sample proportion.

**THEOREM 12:** For all values of $0 < w < 1$, $0 < \alpha < 1$ and $0 \leq \bar{x} \leq 1$ we have:

$$\chi_\alpha^2 \cdot \frac{1 - w}{w} \leq \hat{n}(w, \alpha, \bar{x}) \leq \chi_\alpha^2 \cdot \frac{½ - |w^2 - ½|}{w^2}.$$

It is useful to note that the required sample size function obeys some basic intuitive properties. In Theorem 13 below we show that the function is decreasing in $w$ and $\alpha$. This means that, *ceteris paribus*, obtaining an interval with lower accuracy (higher width) or lower confidence, requires less data, just as we would expect. In Theorem 14 we see that these properties carry over (aside from an exception of weak monotonicity when $w \leq 1/\sqrt{2}$) for the conservative sample size function.

**THEOREM 13:** For all values of $0 < w < 1$, $0 < \alpha < 1$ and $0 \leq \bar{x} \leq 1$ we have:

$$\frac{\partial \hat{n}}{\partial w}(w, \alpha, \bar{x}) < 0 \qquad \frac{\partial \hat{n}}{\partial \alpha}(w, \alpha, \bar{x}) < 0.$$

**THEOREM 14:** Consider the conservative required sample size function:

$$\hat{n}_{\text{cons}}(w, \alpha) \equiv \chi_\alpha^2 \cdot \frac{\frac{1}{2} - |w^2 - \frac{1}{2}|}{w^2}.$$

For all values of $0 < w < 1$ and $0 < \alpha < 1$ we have:

$$\frac{\partial \hat{n}_{\text{cons}}}{\partial w}(w, \alpha) \leq 0 \qquad \frac{\partial \hat{n}_{\text{cons}}}{\partial \alpha}(w, \alpha) < 0.$$

(The first of these monotonicity results is strict for $w > 1/\sqrt{2}$. Note that $\hat{n}_{\text{cons}}$ has a point of discontinuity in its derivative at $w = 1/\sqrt{2}$ but its weak derivative is zero at this point.)

Statistical analysts should give serious thought to the approach they wish to take to sample size calculations for binary data. The conservative required sample size function will give a size that is sufficient for an inference with a prespecified level of accuracy irrespective of the sample data. However, this function is highly conservative and will usually give more accuracy than was necessary. In some statistical problems the analyst will have access to prior information that makes it reasonable to specify an expected sample proportion *a priori* and use the standard (non-conservative) required sample size function. Alternatively, one could work between these extremes, computing the required sample size for a range of expected sample proportions that is smaller than its full range but larger than a single point. In any case, statistical users should reflect on the prior information available in their own work and make a judgement accordingly.

## 5. Summary and conclusion

In this paper we have examined the properties of the Wilson score interval and generalised this confidence interval to inferences for the population proportion or unsampled proportion for the general case where the population can be finite or infinite. Our analysis shows that the "finite population correction" for the latter intervals involves a simple change in the standard formula, to use an "effective sample size" that is fully determined by the sample and population sizes. For both the standard and generalised confidence intervals, we have established a number of intuitive statistics properties for the bounds of the intervals and the widths of the intervals.

Accuracy of the confidence interval can be measured by considering the width of the interval as a function of the sample size and population size, holding the stipulated sample proportion constant. Using this method, for an inference pertaining to the population proportion, we have shown that the width function is a consistent interval estimator, which becomes more accurate as the sample size increases and less accurate as the population size increases. (We consider the accuracy only by the width of the interval, holding the sample proportion as a fixed value; it is of course possible that gaining more data can lead to new data that is not reflective of the general population proportion and so in practice an interval may include the true population proportion for a smaller data set but exclude if for a larger data set.)

For the case where we make an inference about the unsampled proportion, the accuracy of the interval is bounded from below based on the population size, and for a finite population the interval estimator is not consistent. The properties of the interval in this case are more complex, owing to the trade-off between the information gain from more data, and the higher variability in the unsampled part of the population from removing data from that part.

Finally, we have derived the form of the required sample size needed to obtain an interval with a given level of accuracy for a stipulated sample proportion. Varying the sample proportion over its possible range of values gives a conservative form for the required sample size that can be used to conduct an *a priori* sample size calculation to achieve a desired level of accuracy regardless of the sample data. We have seen that the required sample size for inference for the population mean likewise obeys a number of intuitive properties.

The generalised versions of the Wilson score interval considered here are available for use in the `CONF.prop` function in the `stat.extend` package in `R`. This gives a simple and user-friendly facility for statistical practitioners to implement these confidence intervals in contexts where they wish to use binary data to make an inference about a finite population proportion or a finite unsampled proportion. We hope that this paper has been useful to readers in setting out details on the construction and properties of these intervals.

# Appendix I: Proof of Theorems

In this appendix we prove the theorems shown in the main body of the paper.

**PROOF OF THEOREM 1:** As $\alpha \to 0$ we get $\chi_\alpha \to \infty$ so that:

$$\lim_{\alpha \to 0} \mathrm{CI}_\infty(\alpha, n, \bar{x}) = \lim_{\chi_\alpha \to \infty} \left[ \frac{n\bar{x} + \tfrac{1}{2}\chi_\alpha^2}{n + \chi_\alpha^2} \pm \frac{\chi_\alpha}{n + \chi_\alpha^2} \sqrt{n\bar{x}(1-\bar{x}) + \tfrac{1}{4}\chi_\alpha^2} \right]$$

$$= \left[ \lim_{\chi_\alpha \to \infty} \frac{n\bar{x}_n + \tfrac{1}{2}\chi_\alpha^2}{n + \chi_\alpha^2} \pm \lim_{\chi_\alpha \to \infty} \frac{\chi_\alpha}{n + \chi_\alpha^2} \sqrt{n\bar{x}(1-\bar{x}) + \tfrac{1}{4}\chi_\alpha^2} \right]$$

$$= \left[ \lim_{\chi_\alpha \to \infty} \frac{\tfrac{1}{2}\chi_\alpha^2}{\chi_\alpha^2} \pm \lim_{\chi_\alpha \to \infty} \frac{\tfrac{1}{2}\chi_\alpha^2}{\chi_\alpha^2} \right]$$

$$= [\tfrac{1}{2} \pm \tfrac{1}{2}] = [0,1],$$

and as $\alpha \to 1$ we get $\chi_\alpha \to 0$ so that:

$$\lim_{\alpha \to 1} \mathrm{CI}_\infty(\alpha, n, \bar{x}) = \lim_{\chi_\alpha \to 0} \left[ \frac{n\bar{x} + \tfrac{1}{2}\chi_\alpha^2}{n + \chi_\alpha^2} \pm \frac{\chi_\alpha}{n + \chi_\alpha^2} \sqrt{n\bar{x}(1-\bar{x}) + \tfrac{1}{4}\chi_\alpha^2} \right]$$

$$= \left[ \lim_{\chi_\alpha \to 0} \frac{n\bar{x} + \tfrac{1}{2}\chi_\alpha^2}{n + \chi_\alpha^2} \pm \lim_{\chi_\alpha \to 0} \frac{\chi_\alpha}{n + \chi_\alpha^2} \sqrt{n\bar{x}(1-\bar{x}) + \tfrac{1}{4}\chi_\alpha^2} \right]$$

$$= \left[ \frac{n\bar{x}}{n} \pm \frac{0}{n} \sqrt{n\bar{x}(1-\bar{x})} \right]$$

$$= [\bar{x}].$$

As $n \to 0$ we have:

$$\lim_{n \to 0} \mathrm{CI}_\infty(\alpha, n, \bar{x}) = \lim_{n \to 0} \left[ \frac{n\bar{x} + \tfrac{1}{2}\chi_\alpha^2}{n + \chi_\alpha^2} \pm \frac{\chi_\alpha}{n + \chi_\alpha^2} \sqrt{n\bar{x}(1-\bar{x}) + \tfrac{1}{4}\chi_\alpha^2} \right]$$

$$= \left[ \lim_{n \to 0} \frac{n\bar{x} + \tfrac{1}{2}\chi_\alpha^2}{n + \chi_\alpha^2} \pm \lim_{n \to 0} \frac{\chi_\alpha}{n + \chi_\alpha^2} \sqrt{n\bar{x}(1-\bar{x}) + \tfrac{1}{4}\chi_\alpha^2} \right]$$

$$= \left[ \frac{\tfrac{1}{2}\chi_\alpha^2}{\chi_\alpha^2} \pm \frac{\tfrac{1}{2}\chi_\alpha^2}{\chi_\alpha^2} \right]$$

$$= [\tfrac{1}{2} \pm \tfrac{1}{2}] = [0,1],$$

and as $n \to \infty$ we have:

$$\lim_{n \to \infty} \mathrm{CI}_\infty(\alpha, n, \bar{x}) = \lim_{n \to \infty} \left[ \frac{n\bar{x} + \tfrac{1}{2}\chi_\alpha^2}{n + \chi_\alpha^2} \pm \frac{\chi_\alpha}{n + \chi_\alpha^2} \sqrt{n\bar{x}(1-\bar{x}) + \tfrac{1}{4}\chi_\alpha^2} \right]$$

$$= \left[ \lim_{n \to \infty} \frac{n\bar{x} + \tfrac{1}{2}\chi_\alpha^2}{n + \chi_\alpha^2} \pm \lim_{n \to \infty} \frac{\chi_\alpha}{n + \chi_\alpha^2} \sqrt{n\bar{x}(1-\bar{x}) + \tfrac{1}{4}\chi_\alpha^2} \right]$$

$$= \left[ \lim_{n \to \infty} \frac{n\bar{x}}{n} \pm \lim_{n \to \infty} \frac{\chi_\alpha}{n + \chi_\alpha^2} \sqrt{n\bar{x}(1-\bar{x})} \right]$$

$$= \left[ \lim_{n \to \infty} \bar{x} \pm 0 \right]$$

$$= [\bar{x}].$$

This establishes each of the limits in the theorem. ∎

**PROOF OF THEOREM 2:** To facilitate analysis of the first result we will use the derivatives:

$$\phi_\alpha \equiv \frac{\partial \chi_\alpha}{\partial \alpha} \qquad 2\phi_\alpha \chi_\alpha = \frac{\partial \chi_\alpha^2}{\partial \alpha},$$

and we note that $\phi_\alpha < 0$ for all $0 < \alpha < 1$. To obtain the first monotonicity result we have:

$$\frac{\partial w_\infty}{\partial \alpha}(\alpha, n, \bar{x}) = \frac{\partial \chi_\alpha}{\partial \alpha} \cdot \frac{\partial w_\infty}{\partial \chi_\alpha}(\alpha, n, \bar{x})$$

$$= \phi_\alpha \cdot \frac{\partial}{\partial \chi_\alpha} \frac{2\chi_\alpha}{n + \chi_\alpha^2} \sqrt{n\bar{x}(1-\bar{x}) + \tfrac{1}{4}\chi_\alpha^2}$$

$$= \phi_\alpha \left[ \frac{2(n - \chi_\alpha^2)}{(n + \chi_\alpha^2)^2} \sqrt{n\bar{x}(1-\bar{x}) + \tfrac{1}{4}\chi_\alpha^2} + \frac{\chi_\alpha}{n + \chi_\alpha^2} \frac{\tfrac{1}{2}\chi_\alpha}{\sqrt{n\bar{x}(1-\bar{x}) + \tfrac{1}{4}\chi_\alpha^2}} \right]$$

$$= \frac{\phi_\alpha}{\sqrt{n\bar{x}(1-\bar{x}) + \tfrac{1}{4}\chi_\alpha^2}} \left[ \frac{2(n - \chi_\alpha^2)}{(n + \chi_\alpha^2)^2} (n\bar{x}(1-\bar{x}) + \tfrac{1}{4}\chi_\alpha^2) + \frac{\tfrac{1}{2}\chi_\alpha^2}{n + \chi_\alpha^2} \right]$$

$$= \frac{\phi_\alpha}{\sqrt{n\bar{x}(1-\bar{x}) + \tfrac{1}{4}\chi_\alpha^2}} \cdot \frac{2(n - \chi_\alpha^2)(n\bar{x}(1-\bar{x}) + \tfrac{1}{4}\chi_\alpha^2) + \tfrac{1}{2}\chi_\alpha^2(n + \chi_\alpha^2)}{(n + \chi_\alpha^2)^2}$$

$$= \frac{\phi_\alpha}{\sqrt{n\bar{x}(1-\bar{x}) + \tfrac{1}{4}\chi_\alpha^2}} \cdot \frac{2n\bar{x}(1-\bar{x})(n - \chi_\alpha^2) + n\chi_\alpha^2}{(n + \chi_\alpha^2)^2}.$$

We can simplify the numerator in the latter fraction as follows:

$$\begin{aligned}
\text{NUM} &\equiv 2n\bar{x}(1-\bar{x})(n - \chi_\alpha^2) + n\chi_\alpha^2 \\
&= 2n^2 \bar{x}(1-\bar{x}) - 2n\bar{x}(1-\bar{x}_n)\chi_\alpha^2 + n\chi_\alpha^2 \\
&= 2n^2 \bar{x}(1-\bar{x}) + n\chi_\alpha^2(1 - 2\bar{x}(1-\bar{x}_n)) \\
&= 2n^2 \bar{x}(1-\bar{x}) + n\chi_\alpha^2(1 - 2\bar{x} + 2\bar{x}^2) \\
&= 2n^2 \bar{x}(1-\bar{x}) + n\chi_\alpha^2(1 - 2\bar{x} + \bar{x}^2) + n\chi_\alpha^2 \bar{x}^2 \\
&= 2n^2 \bar{x}(1-\bar{x}) + n\chi_\alpha^2(1 - \bar{x})^2 + n\chi_\alpha^2 \bar{x}^2 \\
&= 2n^2 \bar{x}(1-\bar{x}) + n\chi_\alpha^2[\bar{x}^2 + (1-\bar{x})^2] \\
&= 2n^2 \bar{x}(1-\bar{x}) + n\chi_\alpha^2[\bar{x}^2 + (1-\bar{x})^2].
\end{aligned}$$

We therefore obtain:

$$\frac{\partial w_\infty}{\partial \alpha}(\alpha, n, \bar{x}) = \phi_\alpha \cdot \frac{2n^2 \bar{x}(1-\bar{x}) + n\chi_\alpha^2[\bar{x}^2 + (1-\bar{x})^2]}{(n + \chi_\alpha^2)^2 \sqrt{n\bar{x}(1-\bar{x}) + \tfrac{1}{4}\chi_\alpha^2}} < 0.$$

(The final inequality is obtained by observing that all the terms in the fractional part are positive and the value $\phi_\alpha$ is negative.) To obtain the second monotonicity result we have:

$$\frac{\partial w_\infty}{\partial n}(\alpha, n, \bar{x}) = \frac{\partial}{\partial n} \frac{2\chi_\alpha}{n + \chi_\alpha^2} \sqrt{n\bar{x}(1-\bar{x}) + \tfrac{1}{4}\chi_\alpha^2}$$

$$= -\frac{2\chi_\alpha}{(n+\chi_\alpha^2)^2} \sqrt{n\bar{x}(1-\bar{x}) + \tfrac{1}{4}\chi_\alpha^2} + \frac{\chi_\alpha}{n+\chi_\alpha^2} \frac{\bar{x}(1-\bar{x})}{\sqrt{n\bar{x}(1-\bar{x}) + \tfrac{1}{4}\chi_\alpha^2}}$$

$$= -\frac{\chi_\alpha}{(n+\chi_\alpha^2)^2} \cdot \frac{2(n\bar{x}(1-\bar{x}) + \tfrac{1}{4}\chi_\alpha^2) - (n+\chi_\alpha^2)\bar{x}(1-\bar{x})}{\sqrt{n\bar{x}(1-\bar{x}) + \tfrac{1}{4}\chi_\alpha^2}}$$

$$= -\frac{\chi_\alpha}{(n+\chi_\alpha^2)^2} \cdot \frac{\tfrac{1}{2}\chi_\alpha^2 + n\bar{x}(1-\bar{x}) - \bar{x}(1-\bar{x})\chi_\alpha^2}{\sqrt{n\bar{x}(1-\bar{x}) + \tfrac{1}{4}\chi_\alpha^2}}$$

$$= -\frac{\chi_\alpha}{(n+\chi_\alpha^2)^2} \cdot \frac{n\bar{x}(1-\bar{x}) + \tfrac{1}{2}\chi_\alpha^2(1 - 2\bar{x} + 2\bar{x}^2)}{\sqrt{n\bar{x}(1-\bar{x}) + \tfrac{1}{4}\chi_\alpha^2}}$$

$$= -\frac{\chi_\alpha}{(n+\chi_\alpha^2)^2} \cdot \frac{n\bar{x}(1-\bar{x}) + \tfrac{1}{2}\chi_\alpha^2[\bar{x}^2 + (1-\bar{x})^2]}{\sqrt{n\bar{x}(1-\bar{x}) + \tfrac{1}{4}\chi_\alpha^2}} < 0.$$

(The final inequality is obtained by observing that all the terms in both the fractional parts are positive so the negative sign at the front determines the sign.) ∎

**PROOF OF THEOREM 3:** We have:

$$\frac{\partial L_\infty}{\partial \bar{x}}(\alpha, n, \bar{x}) = \frac{\partial}{\partial \bar{x}} \frac{n\bar{x} + \tfrac{1}{2}\chi_\alpha^2}{n + \chi_\alpha^2} - \frac{\partial}{\partial \bar{x}} \frac{\chi_\alpha}{n + \chi_\alpha^2} \sqrt{n\bar{x}(1-\bar{x}) + \tfrac{1}{4}\chi_\alpha^2}$$

$$= \frac{n}{n + \chi_\alpha^2} - \frac{\chi_\alpha}{n + \chi_\alpha^2} \frac{n(\tfrac{1}{2} - \bar{x})}{\sqrt{n\bar{x}(1-\bar{x}) + \tfrac{1}{4}\chi_\alpha^2}}$$

$$= \frac{n}{n + \chi_\alpha^2} \cdot \frac{\sqrt{n\bar{x}(1-\bar{x}) + \tfrac{1}{4}\chi_\alpha^2} - (\tfrac{1}{2} - \bar{x})\chi_\alpha}{\sqrt{n\bar{x}(1-\bar{x}) + \tfrac{1}{4}\chi_\alpha^2}}$$

$$= \frac{n}{n + \chi_\alpha^2} \cdot \frac{\sqrt{n\bar{x}(1-\bar{x}) + \tfrac{1}{4}\chi_\alpha^2} - (\tfrac{1}{2} - \bar{x})\chi_\alpha}{\sqrt{n\bar{x}(1-\bar{x}) + \tfrac{1}{4}\chi_\alpha^2}}.$$

To establish that the numerator in the latter fraction is positive, we first note that:

$$n\bar{x}(1-\bar{x}) + \tfrac{1}{4}\chi_\alpha^2 - (\tfrac{1}{2} - \bar{x})^2 \chi_\alpha = n\bar{x}(1-\bar{x}) + \tfrac{1}{4}\chi_\alpha^2 - (\tfrac{1}{4} - \bar{x} + \bar{x}^2)\chi_\alpha$$

$$= n\bar{x}(1-\bar{x}) - (-\bar{x} + \bar{x}^2)\chi_\alpha$$

$$= n\bar{x}(1-\bar{x}) + \bar{x}(1-\bar{x}^2)\chi_\alpha$$

$$= (n + \chi_\alpha)\bar{x}(1-\bar{x}) > 0.$$

This implies that:

$$n\bar{x}(1-\bar{x}) + \tfrac{1}{4}\chi_\alpha^2 > (\tfrac{1}{2} - \bar{x})^2 \chi_\alpha,$$

which implies that:

$$\sqrt{n\bar{x}(1-\bar{x}) + \tfrac{1}{4}\chi_\alpha^2} > (\tfrac{1}{2} - \bar{x})\chi_\alpha.$$

From this latter inequality we then have:

$$\frac{\partial L_\infty}{\partial \bar{x}}(\alpha, n, \bar{x}) = \frac{n}{n + \chi_\alpha^2} \cdot \frac{\sqrt{n\bar{x}(1-\bar{x}) + \tfrac{1}{4}\chi_\alpha^2} - (\tfrac{1}{2} - \bar{x})\chi_\alpha}{\sqrt{n\bar{x}(1-\bar{x}) + \tfrac{1}{4}\chi_\alpha^2}} > 0.$$

The corresponding proof for monotonicity of the upper bound function is virtually identical, and we leave it to the interested reader to establish. To obtain the concavity result we have:

$$\frac{\partial w_\infty}{\partial \bar{x}}(\alpha, n, \bar{x}) = \frac{2\chi_\alpha}{n + \chi_\alpha^2} \frac{n(\tfrac{1}{2} - \bar{x})}{\sqrt{n\bar{x}(1-\bar{x}) + \tfrac{1}{4}\chi_\alpha^2}},$$

$$\frac{\partial^2 w_\infty}{\partial \bar{x}^2}(\alpha, n, \bar{x}) = -\frac{2\chi_\alpha}{n + \chi_\alpha^2} \cdot \frac{n(n\bar{x}(1-\bar{x}) + \tfrac{1}{4}\chi_\alpha^2) + n^2(\tfrac{1}{2} - \bar{x})^2}{(n\bar{x}(1-\bar{x}) + \tfrac{1}{4}\chi_\alpha^2)^{3/2}}$$

$$= -\frac{2\chi_\alpha}{n + \chi_\alpha^2} \cdot \frac{n^2[\bar{x}(1-\bar{x}) + (\tfrac{1}{2} - \bar{x})^2] + \tfrac{1}{4}n\chi_\alpha^2}{(n\bar{x}(1-\bar{x}) + \tfrac{1}{4}\chi_\alpha^2)^{3/2}}$$

$$= -\frac{2\chi_\alpha}{n + \chi_\alpha^2} \cdot \frac{n^2[(\bar{x} - \bar{x}^2) + (\tfrac{1}{4} - \bar{x} + \bar{x}^2)] + \tfrac{1}{4}n\chi_\alpha^2}{(n\bar{x}(1-\bar{x}) + \tfrac{1}{4}\chi_\alpha^2)^{3/2}}$$

$$= -\frac{\chi_\alpha}{n + \chi_\alpha^2} \cdot \frac{\tfrac{1}{2}n(n + \chi_\alpha^2)}{(n\bar{x}(1-\bar{x}) + \tfrac{1}{4}\chi_\alpha^2)^{3/2}} < 0.$$

This establishes each of the results in the theorem. ∎

**PROOF OF THEOREM 4:** It is trivial to establish that $w_\infty(\alpha, n, \bar{x}) \geq 0$, which establishes the inequality $L_\infty(\alpha, n, \bar{x}) \leq U_\infty(\alpha, n, \bar{x})$. The remaining bounds in the theorem can be obtained using the fact that that both interval bounds are monotonically increasing in $\bar{x}$ (Theorem 3). This monotonicity property means that the minimum value of the lower bound function is achieved when $\bar{x} = 0$, giving:

$$L_\infty(\alpha, n, \bar{x}) \geq L_\infty(\alpha, n, 0) = \frac{\tfrac{1}{2}\chi_\alpha^2}{n + \chi_\alpha^2} - \frac{\tfrac{1}{2}\chi_\alpha^2}{n + \chi_\alpha^2} = 0,$$

and the maximum value of the upper bound function is achieved when $\bar{x} = 1$, which gives:

$$U_\infty(\alpha, n, \bar{x}) \leq U_\infty(\alpha, n, 1) = \frac{n + \tfrac{1}{2}\chi_\alpha^2}{n + \chi_\alpha^2} + \frac{\tfrac{1}{2}\chi_\alpha^2}{n + \chi_\alpha^2} = 1.$$

Putting these results together establishes the bounds shown in the theorem. ∎

**PROOF OF THEOREM 5:** Theorem 3 establishes that the width function is strictly concave. Since the width function depends on $\bar{x}$ only through $\bar{x}(1-\bar{x})$, it is symmetric in $\bar{x}$ around the mid-point $\bar{x} = ½$. Consequently, the function is minimised at its extremes and maximised at its mid-point; substitution of these values gives the bounds shown in the theorem. ∎

**PROOF OF THEOREM 6:** The special values of the effective sample size functions follow easily by substitution, and can be confirmed by the reader. For the monotonicity/directional results we have:

$$\frac{\partial n_*}{\partial n}(n, N) = \frac{N(N-1)}{(N-n)^2} > 0,$$

$$\frac{\partial n_*}{\partial N}(n, N) = -\frac{n(n-1)}{(N-n)^2} < 0,$$

$$\frac{\partial n_{**}}{\partial n}(n, N) = \frac{N - 2n}{N - 1},$$

$$\frac{\partial n_{**}}{\partial N}(n, N) = \frac{n-1}{(N-1)^2} > 0,$$

and the sign for the third derivative follows trivially. ∎

**PROOF OF THEOREM 7:** These results follow easily by combining Theorems 1 and 6. ∎

**PROOF OF THEOREM 8:** The results in the first column follow easily from Theorem 2. For the remaining results we apply the chain rule and apply the monotonicity results in Theorems 2 and 6. For the width function $w_N$ we obtain:

$$\frac{\partial w_N}{\partial n}(\alpha, n, \bar{x}) = \frac{\partial}{\partial n} w_\infty(\alpha, n_*(n, N), \bar{x})$$

$$= \frac{\partial n_*}{\partial n}(n, N) \cdot \frac{\partial w_\infty}{\partial n}(\alpha, n_*, \bar{x}) < 0,$$

$$\frac{\partial w_N}{\partial N}(\alpha, n, \bar{x}) = \frac{\partial}{\partial N} w_\infty(\alpha, n_*(n, N), \bar{x})$$

$$= \frac{\partial n_*}{\partial N}(n, N) \cdot \frac{\partial w_\infty}{\partial n}(\alpha, n_*, \bar{x}) > 0.$$

For the width function $w_{n:N}$ we obtain:

$$\frac{\partial w_{n:N}}{\partial n}(\alpha, n, \bar{x}) = \frac{\partial}{\partial n} w_\infty(\alpha, n_{**}(n, N), \bar{x})$$

$$= \frac{\partial n_{**}}{\partial n}(n, N) \cdot \frac{\partial w_\infty}{\partial n}(\alpha, n_*, \bar{x}),$$

$$\frac{\partial w_{n:N}}{\partial N}(\alpha, n, \bar{x}) = \frac{\partial}{\partial N} w_\infty(\alpha, n_{**}(n, N), \bar{x})$$

$$= \frac{\partial n_{**}}{\partial N}(n, N) \cdot \frac{\partial w_\infty}{\partial n}(\alpha, n_*, \bar{x}) < 0.$$

and the sign result is obtained as:

$$\operatorname{sgn}\frac{\partial w_{n:N}}{\partial n} = \operatorname{sgn}\frac{\partial n_{**}}{\partial N}(n, N) \times \operatorname{sgn}\frac{\partial w_\infty}{\partial n}(\alpha, n_*, \bar{x})$$
$$= \operatorname{sgn}(\tfrac{1}{2}N - n) \times (-1)$$
$$= \operatorname{sgn}(n - \tfrac{1}{2}N),$$

which completes the proof. ■

**PROOF OF THEOREM 9:** These results follow easily from Theorem 3. ■

**PROOF OF THEOREM 10:** These results follow easily from Theorem 4. ■

**PROOF OF THEOREM 11:** The upper bound is easily established from theorem 7, using the fact that $\operatorname{CI}_N(\alpha, 0, \bar{x}) = [0,1]$. For the lower bound, we first note that $n_{**}$ is a strictly concave function in $n$ and is maximised at the critical point $\hat{n} = \tfrac{1}{2}N$, which gives the maximum value $\max_{0 \leq n \leq N} n_{**}(n, N) = \tfrac{1}{4}N$. Since $w_{n:N}(\alpha, n, \bar{x}) = w_\infty(\alpha, n_{**}, \bar{x})$, application of Theorem 2 shows that the width function for the confidence interval is monotonically decreasing in $n_{**}$, so this choice of sample size gives the minimum of the width function:

$$\min_{0 \leq n \leq N} w_{n:N}(\alpha, n, \bar{x}) = w_\infty(\alpha, \tfrac{1}{4}N, \bar{x}) = \frac{\chi_\alpha}{\tfrac{1}{4}N + \chi_\alpha^2}\sqrt{N\bar{x}(1 - \bar{x}) + \chi_\alpha^2}.$$

This function is minimised by taking $\bar{x} = \tfrac{1}{2}$, which gives the absolute minimum:

$$\min_{0 \leq n \leq N} \min_{0 \leq \bar{x} \leq 1} w_{n:N}(\alpha, n, \bar{x}) = \frac{\chi_\alpha^2}{\tfrac{1}{4}N + \chi_\alpha^2}.$$

This establishes the lower bound in the theorem, which completes the proof. ■

**PROOF OF THEOREM 12:** Letting $z \equiv \bar{x}(1 - \bar{x})$ we have:

$$\hat{n}(w, \alpha, \bar{x}) = \frac{\chi_\alpha^2}{w^2}\left[\sqrt{w^2 - 4zw^2 + 4z^2} - \sqrt{w^4 - 4zw^2 + 4z^2}\right].$$

Differentiating with respect to $z$ we obtain:

$$\frac{\partial \hat{n}}{\partial z}(w, \alpha, \bar{x}) = \frac{2\chi_\alpha^2}{w^2}\left[\frac{w^2 - 2z}{\sqrt{w^4 - 4zw^2 + 4z^2}} - \frac{w^2 - 2z}{\sqrt{w^2 - 4zw^2 + 4z^2}}\right] > 0.$$

Since $\hat{n}$ is increasing in $z$ it follows that the required sample size is minimised at the values $\bar{x} = 0$ and $\bar{x} = 1$ (corresponding to $z = 0$) and maximised at the value $\bar{x} = \tfrac{1}{2}$ (corresponding to $z = \tfrac{1}{4}$). Substitution of these values yields the bounds in the theorem. ■

**PROOF OF THEOREM 13:** It is possible to prove the first result with standard calculus methods, but here we will prove it more directly. Since $w_\infty$ is strictly decreasing in $n$, for any $w < w'$ we have:

$$\hat{n}(w, \alpha, \bar{x}) = \inf\{n \geq 0 | w_\infty(\alpha, n, \bar{x}) \leq w\}$$
$$< \inf\{n \geq 0 | w_\infty(\alpha, n, \bar{x}) \leq w'\} = \hat{n}(w', \alpha, \bar{x}).$$

Since $\hat{n}$ is differentiable, the first result in the theorem follows. To prove the second result, we first note that $\partial \chi_\alpha^2 / \partial \alpha < 0$. Applying the chain rule, we therefore have:

$$\frac{\partial \hat{n}}{\partial \alpha}(w, \alpha, \bar{x}) = \frac{\partial \chi_\alpha^2}{\partial \alpha} \cdot \frac{\partial \hat{n}}{\partial \chi_\alpha^2}(w, \alpha, \bar{x}) = \frac{\partial \chi_\alpha^2}{\partial \alpha} \cdot \frac{\hat{n}(w, \alpha, \bar{x})}{\chi_\alpha^2} < 0,$$

which establishes the theorem. ∎

**PROOF OF THEOREM 14:** Since $\chi_\alpha^2$ is decreasing in $\alpha$ the second monotonicity property in the theorem is trivial to establish. To establish the first property, we first note that we can rewrite the conservative sample size function as:

$$\hat{n}_{\text{cons}}(w, \alpha) = \begin{cases} \chi_\alpha^2 & \text{for } w \leq \frac{1}{\sqrt{2}}, \\ \chi_\alpha^2 \cdot \frac{1-w^2}{w^2} & \text{for } w > \frac{1}{\sqrt{2}}. \end{cases}$$

For $w \leq 1/\sqrt{2}$ the function is constant with respect to $w$. For $w > 1/\sqrt{2}$ we have:

$$\frac{\partial \hat{n}_{\text{cons}}}{\partial w}(w, \alpha) = \chi_\alpha^2 \cdot \frac{\partial}{\partial w} \frac{1-w^2}{w^2}$$
$$= \chi_\alpha^2 \cdot \frac{-2w^3 - 2(1-w^2)w}{w^4}$$
$$= -\frac{2\chi_\alpha^2}{w^3} < 0.$$

Putting these properties together establishes the result in the theorem. ∎

# Appendix II: Alternative form for the intervals

**ALTERNATIVE FORM:** We can state the confidence intervals is in terms of the parameters:

$$\phi_* \equiv \frac{N-n}{N} \cdot \frac{\chi_\alpha^2}{2n} \qquad \phi_{**} \equiv \frac{N}{N-n} \cdot \frac{\chi_\alpha^2}{2n}.$$

Using these parameters we have:

$$n_* = \frac{\chi_\alpha^2}{2\phi_*} \qquad n_{**} = \frac{\chi_\alpha^2}{2\phi_*}.$$

With this alternative parameterisation we can write the confidence intervals as:

$$\mathrm{CI}_N(\alpha, n, \bar{x}) = \left[ \frac{\bar{x} + \phi_*}{1 + 2\phi_*} \pm \frac{1}{1 + 2\phi_*} \sqrt{2\phi_* \bar{x}(1 - \bar{x}) + \phi_*^2} \right],$$

$$\mathrm{CI}_{n:N}(\alpha, n, \bar{x}) = \left[ \frac{\bar{x} + \phi_{**}}{1 + 2\phi_{**}} \pm \frac{1}{1 + 2\phi_{**}} \sqrt{2\phi_{**} \bar{x}(1 - \bar{x}) + \phi_{**}^2} \right].$$

This parameterisation collects the parameters $n$, $N$ and $\alpha$ into a single parameter $\phi_*$ or $\phi_{**}$ and then gives the confidence interval in a succinct form that separates $\bar{x}$ from the other parameters. This parameterisation is used in constructing the **CONF.prop** function in the **stat.extend** package in **R** (O'Neill and Fultz 2020)